\documentclass[reqno]{amsart}

\usepackage{amssymb,amsmath}
\usepackage{euscript}
\usepackage{graphicx}
\usepackage{tikz}
\usepackage{tikz-cd}
\usepackage{hyperref}
\usepackage{mathrsfs}
\usepackage{dsfont}
\usepackage{mathtools}
\usepackage{microtype}
\usepackage{booktabs}
\usepackage{array}

\newtheorem{thm}{Theorem}[section]
\newtheorem{corollary}[thm]{Corollary}
\newtheorem{lemma}[thm]{Lemma}
\newtheorem{proposition}[thm]{Proposition}
\usepackage{adjustbox}
\usepackage{float}

\numberwithin{equation}{section}
\renewcommand{\subsection}[1]{\hspace{-\parindent}\refstepcounter{subsection}{\bf (\thesubsection) #1.}}
\renewcommand{\thesubsection}{\arabic{section}\alph{subsection}}

\allowdisplaybreaks

\newtheorem{def-lem}[thm]{Definition-Lemma}

\theoremstyle{definition}
\newtheorem{definition}[thm]{Definition}

\theoremstyle{remark}

\newtheorem{remark}[thm]{Remark}

\numberwithin{equation}{section}

\title{On the strong Arnold conjecture for Grassmannians}
\author{Kenneth Blakey}
\address{Department of Mathematics, MIT, 182 Memorial Drive, Cambridge, MA 02139, U.S.A.} \email{kblakey@mit.edu, romakrut@mit.edu}

\author{Roman Krutowski}

\begin{document}

\begin{abstract}
On several families of flag varieties, we provide examples of degenerate Hamiltonian diffeomorphisms exhibiting significantly smaller numbers of fixed points than predicted by the strong Arnold conjecture.
\end{abstract}

\maketitle

\begin{quote}
\emph{Mistakes are an important and instructive part of mathematics, perhaps as important a part as proofs. Proofs are to mathematics what spelling (or even calligraphy) is to poetry. Mathematical works consist of proofs, just as poems consist of characters.} -- \emph{Vladimir Igorevich Arnold}
\end{quote}

\section{Introduction}
\subsection{The Arnold number}
The strong Arnold conjecture, initially stated in \cite[Remark C]{Arnold1965surune} or \cite[Appendix 10]{Arnold1978methods}, predicted that the minimal number of fixed points of a Hamiltonian diffeomorphism on a closed symplectic manifold $(M, \omega)$, which we refer to as the \emph{Arnold number} $\operatorname{Arn}(M)$, is at least the minimal number of critical points of a smooth function on $M$:
  \begin{equation}
  \operatorname{Arn}(M) \ge \mathrm{Crit}(M).
  \end{equation}
  
We have learned from recent works that this lower bound does not hold in general, cf. \cite{jiao2026twofixed, OAIdeg26}. Nevertheless, the question of determining $\operatorname{Arn}(M)$ for an arbitrary $M$ is still completely open. In this work, we provide several infinite families satisfying
  \begin{equation}
  \mathrm{Arn}(M)\leq\mathrm{Crit}(M)-N
  \end{equation}
for any given $N>0$, as well as an infinite family of examples satisfying 
  \begin{equation}
  \operatorname{Arn}(M) \lesssim \frac{\mathrm{Crit}(M)}{2}.
  \end{equation}
  The more precise statements are given in Theorem~\ref{thm: main theorem}.
The symplectic manifolds in question all belong to the class of complex flag varieties, a class whose topology, symplectic geometry, and Floer theory have been quite well studied -- this goes back to the work of Bott and Samelson \cite{BottSamelson1958}. Our expectation is that these examples will be instructive for further progress towards determining $\operatorname{Arn}(M)$ for a sufficiently general symplectic manifold $M$.

\subsection{Some history}
Unfortunately, the following is most likely non-exhaustive, and any omission is the fault of the authors. 

As one might argue \cite{audin2014}, the Arnold conjecture gave birth to symplectic geometry, and consequently sparked the creation of Floer theory \cite{Floer1986,Flo88b, Flo89a, Flo89b}. Tremendous progress has been made towards showing the non-degenerate version of the Arnold conjecture, and its \emph{homological} version has been resolved for both general symplectic manifolds and over various coefficient rings, cf. \cite{LT98, FO99, Pardon16, AB21, FW22, BX22, Rez22}.

In the degenerate case, much less is known. Fortune \cite{Fortune1985} realized the lower bound for complex projective space $\mathbb{C}P^n$. Meanwhile, for tori $T^{2n}$, Conley and Zehnder produced a proof \cite{CZ1983}. Already for products $T^{2m} \times \mathbb{C}P^k$, Oh \cite{Oh1990symplectic} could only produce a lower bound of $k+1$ fixed points, and suggested a counterexample should already appear for $S^2 \times S^2$ or $T^2 \times S^2$. Eliashberg and Floer resolved the case of Riemann surfaces \cite{Eliashberg1978, Floer1986}, and further progress was made by Rudyak and Oprea \cite{RO1999} who proved the conjecture in the symplectically aspherical case; the latter generalized earlier works by Hofer \cite{Hof88} and Floer \cite{Flo89a}. We should also mention \cite{LeOno1994} as it played a crucial role in the development of the theory. Finally, the strong Arnold conjecture was recently disproved, using generative AI, in \cite{jiao2026twofixed, OAIdeg26} via two different counterexamples that we discuss below.

Note, the strong Arnold conjecture for Hamiltonian diffeomorphisms sufficiently close to the identity is also contained in \cite[Appendix 9]{Arnold1978methods}; as it turns out, this statement is closer to the correct one. Significant progress has been made to prove this form of the conjecture, and we refer the reader to \cite{Weinstein1983, LS1985, chaperon1983, Hof85, AlbersHein2016, AAC2025, AsselleStarostka26} and the references therein. 

We should also mention that any closed symplectic manifold of dimension at least 4 admits a Hamiltonian homeomorphism with a single fixed point, cf. \cite{BHS2018}.

Although the question of determining the Arnold number remains largely open, there has been significant progress toward obtaining general lower bounds for $\mathrm{Arn}(M)$. In the case of symplectic toric manifolds, Givental \cite{Giv95} provides lower bounds for any symplectomorphism (in the monotone case, the lower bound is given by the minimal Chern number). Using quantum cohomology, Schwarz \cite{Schwarz1998QuantumCupLength} obtains lower bounds for a given symplectomorphism in terms of both the Hofer norm and quantum cuplength; in particular, on every positively monotone symplectic manifold the authors know of, this achieves a lower bound given by the minimal Chern number via taking a degree 2 non-nilpotent element in quantum cohomology. Moreover, recent work of Gong \cite{Gong2025DegenerateFixedPoints} improves upon Schwarz's results and provides lower bounds in terms of Ruan-Tian's mixed Gromov-Witten invariants \cite{RuanTian1995QuantumCohomology}. In particular, we would like to emphasize that, for the examples we consider in this paper, the best lower bound that the authors are aware of is given via the minimal Chern number as the quantum cohomology has been computed.

\subsection{Main results}
 We provide many counterexamples to the strong Arnold conjecture, and we summarize our findings as the following result, which follows from Corollaries~\ref{cor: arn for LG}, \ref{cor: arn for type B}, and Proposition~\ref{prop: quadrics}.

\begin{thm}\label{thm: main theorem}
  For a complex flag variety $M$, in Table \ref{tab:arn-bounds}, the Arnold number $\operatorname{Arn}(M)$ satisfies the upper and lower bounds listed in the aforementioned table.\footnote{Here, $OG(n-k, 2n+1)$ is the type $B$ isotropic Grassmannian and $LG(n, 2n)$ is the complex Lagrangian Grassmannian; note, $n, k \in \mathbb{Z}_{\ge 0}, \, 0 \leq k < n$.}
\end{thm}

\begin{remark}
As already mentioned, every lower bound in Table \ref{tab:arn-bounds} is given via the minimal Chern number (in particular, we would like to emphasize that a lower bound given by 4 holds in the case of $OG(1,6)$, even though the table excludes $n=4$). The minimal Chern lower bound follows via \cite[Corollary 1.3]{Schwarz1998QuantumCupLength} together with a computation of the relevant quantum cohomology ring; the latter is computed for each row in \cite{Tamvakis2005IsotropicGrassmannians}.
\end{remark}

  \begin{table}[H]
 \renewcommand{\arraystretch}{1.5}
 \setlength{\tabcolsep}{1ex}
 \begin{adjustbox}{center}
 \begin{tabular}
 {|>{$}l<{$}|>{$}c<{$}|>{$}c<{$}|>{$}c<{$}|}
  \hline
  M & \mathrm{Crit}(M) & \operatorname{Arn}(M) \le & \operatorname{Arn}(M) \ge \\
  \hline\hline
   OG(n-k, 2n+1), \, k>0 & (n-k)(n+k+1)-\binom{n-k+1}{2}+1 & \mathrm{Crit}(M)-(n-k) & n+k \\ \hline
   OG(n, 2n+1)& \binom{n+1}{2}+1 & \lfloor{\tfrac{n^2}{2}\rfloor}+2 & 2n \\ \hline
   OG(1, n+2) \cong Q^n, \, n \neq 4& n+1& 2\; (n=1),\;n & 2\; (n=1),\;n \\ \hline
   LG(n, 2n), \, n\;\mathrm{even}& \binom{n+1}{2}+1 & \tfrac{n(n+2)}{4}+1 & n+1 \\ \hline
   LG(n, 2n), \, n\;\mathrm{odd}& \binom{n+1}{2}+1 & \tfrac{(n+1)^2}{4}+1 & n+1 \\ \hline
 \end{tabular}
 \end{adjustbox}
 \caption{Bounds for $\operatorname{Arn}(M)$ compared with $\mathrm{Crit}(M)$.}
 \label{tab:arn-bounds}
\end{table}
 \subsubsection{Degenerate Hamiltonian diffeomorphisms}
 We briefly comment on the strategy of our proof, which is a slight generalization of the strategy in \cite{OAIdeg26}; it begins with the following observation: given any \emph{clean} Hamiltonian diffeomorphism $\varphi$, one can always slightly perturb it to get a (potentially degenerate) Hamiltonian diffeomorphism $\tilde{\varphi}$ with exactly
  \begin{equation}
  \mathrm{Crit}\big(\mathrm{Fix}(\varphi)\big)
  \end{equation}
many fixed points, where $\mathrm{Fix}(\varphi)$ is the submanifold of fixed points of $\varphi$. A more formal statement is given in Proposition~\ref{proposition: trick}. The version of this statement for finite-order Hamiltonian diffeomorphisms is the main idea in the construction for $Q^3$ in \cite{OAIdeg26}.

 The second ingredient comes from the fact that any complex flag variety, corresponding to a given simple Lie group, comes equipped with a Hamiltonian action of its maximal compact subgroup. In the case of type $B$ isotropic Grassmannians $M=OG(n-k, 2n+1)$, it is the group $G=SO(2n+1)$; in the case of type $C$ complex Lagrangian Grassmannians $M=LG(n,2n)$, it is the compact symplectic group $G=Sp(n)$. These actions are very well understood, and each element $g \in G$ with 
 \begin{equation}
 \mathrm{Crit}(\mathrm{Fix}\big(g)\big) < \mathrm{Crit}(M)= \operatorname{dim}_{\mathbb{C}}M+1
 \end{equation}
 provides a counterexample to the strong Arnold conjecture in its original formulation. Moreover, in Corollaries~\ref{cor: arn for LG} and \ref{cor: arn for type B}, we obtain the optimal bound among all elements in these $G$.

 The optimization among elements of these compact groups is done via an inductive procedure of ``stretching films'' between the fixed point set of a given element $g$ by the fixed point set of a more optimal element $\tilde{g}$. We anticipate that this procedure may have generalizations beyond the settings of Lie group actions on flag varieties, and may play a role in obtaining optimal Hamiltonian diffeomorphisms for other more general symplectic manifolds.

 We should mention that there is an alternative strategy being developed by Abouzaid and Shelukhin which suggests that one should always be able to merge fixed points connected by ``twisted isotropic" spheres, and this is precisely the strategy implemented by \cite{jiao2026twofixed} in the case of $Q^2= \mathbb{P}^1 \times \mathbb{P}^1$. It may very well turn out that these two strategies can complement each other for the general question of determining the Arnold number.

\emph{Acknowledgments}
We would like to thank Paul Seidel for fruitful discussions and helpful suggestions, and Shaoyun Bai for reading an early draft. The first author was partially supported by an NSF Graduate Research Fellowship award during this work.

\emph{AI Disclosure} 
In this project, we used large language models to perform several computations (some of which appear in the proofs of Proposition~\ref{prop: LG Arnold} and Proposition~\ref{prop: arn for type B}), search the literature, and critically review our arguments. The ideas behind Proposition~\ref{proposition: trick}, and generalizing from quadrics to complex flag varieties, arose naturally through conversations between the authors after studying the manuscript \cite{OAIdeg26}. All writing is that of the authors.

\section{Degenerate Hamiltonian diffeomorphisms on flag varieties}
\subsection{Perturbing clean Hamiltonian diffeomorphisms}
\begin{definition}
  We say that a Hamiltonian diffeomorphism $\varphi \in \operatorname{Ham}(M)$ is \emph{clean} if the diagonal $\Delta_M \subset M \times M$ and the graph 
  \begin{equation}
  \Gamma_\varphi=\Big\{\big(x, \varphi(x)\big) \in M \times M : x\in M\Big\}
  \end{equation}
  intersect cleanly. Note, this is equivalent to saying that each connected component $N$ of the fixed point set $\mathrm{Fix}(\varphi)$ is a Morse-Bott submanifold, i.e.,
  \begin{equation}
    T_xN=\operatorname{ker}(d\varphi_x-\operatorname{Id}),\;\;\forall x\in N.
  \end{equation}
\end{definition}

All the constructions we provide stem from the following quite simple observation.
\begin{proposition}\label{proposition: trick}
  Given a clean Hamiltonian diffeomorphism $\varphi \in \operatorname{Ham}(M)$ with fixed point set given by a submanifold $N\subset M$, we have
  \begin{equation}
    \operatorname{Arn}(M) \le \mathrm{Crit}(N).
  \end{equation}
\end{proposition}

\begin{proof}
We show that $\varphi$ can be slightly adjusted so that the resulting Hamiltonian diffeomorphism has exactly $\mathrm{Crit}(N)$ fixed points. For simplicity of the proof, we assume $N$ to be connected; the disconnected case follows immediately as all of our constructions are local with respect to $N$ (in particular, the new Hamiltonian diffeomorphism coincides with $\varphi$ away from a tubular neighborhood of $N$).

\emph{Step 1: Hamiltonian isotopy in Pozniak's model:} By Pozniak's result \cite[Proposition 3.4.1]{pozniak1999}, there is a symplectomorphism $\Psi$ from a tubular neighborhood $U$ of $N \subset T^*M$ to a neighborhood $\Psi(U)$ of the diagonal 
    \begin{equation}
    \Delta_N \subset\big(M \times M, \Omega=\omega \oplus (-\omega)\big)
    \end{equation}
such that
\begin{equation}
  (M|_U, \nu_N|_U) =\big(\Psi^{-1}(\Delta_M), \Psi^{-1}(\Gamma_\varphi)\big),
\end{equation}
where $\nu_N$ denotes the conormal bundle of $N$. Note, this Darboux-type theorem applies to any pair of cleanly intersecting Lagrangians, and we recall that the graph of a clean Hamiltonian diffeomorphism intersects the diagonal cleanly.

Let $G$ be a function on $M|_U$ which: has exactly $\mathrm{Crit}(N)$ critical points in a neighborhood of $N$, all belonging to $N$, and vanishes on the complement of that neighborhood. Additionally, we choose a cutoff function $\chi$ on $T^*M$ with support inside $U$ and satisfying $\chi \equiv 1$ on a neighborhood of $N$. Possibly scaling $G$ by a small constant, we obtain a Hamiltonian isotopy $\Phi_t$ of $T^*M$, given by
\begin{equation}
(q, p) \mapsto (q, p+t dG),
\end{equation}
in a neighborhood of $N$ with Hamiltonian function $H_t=\chi \cdot (G \circ \pi_{T^*M})$. Now: \begin{equation}
\#\big(\Phi_1(\nu_N) \cap M\big)=\mathrm{Crit}(N).
\end{equation}
Let us also note that $H_t^\Psi=H_t \circ \Psi^{-1}$ is a well-defined Hamiltonian function on $M \times M$, and we denote the flow of its Hamiltonian vector field by $\Phi_t^{\Psi}$. 

\emph{Step 2: from Hamiltonian isotopy to Hamiltonian diffeomorphism:} What remains to show is that $\Phi_1^\Psi(\Gamma_\varphi)$ is the graph of a Hamiltonian diffeomorphism. First of all, it is straightforward to see that there is a family of symplectomorphisms $\varphi_t$, with $\varphi_0=\varphi$, such that $\Phi_t^\Psi(\Gamma_\varphi)=\Gamma_{\varphi_t}$. Indeed, we can set
    \begin{equation}
        \varphi_t=\Phi_{t2}\circ \Phi_{t1}^{-1},\;\;\mathrm{where}\;\Phi_{t1}=\pi_1 \circ \Phi_t^{\Psi} \circ j_\varphi,\;\Phi_{t2}=\pi_2 \circ \Phi_t^{\Psi} \circ j_\varphi,
    \end{equation}
    with $\pi_i$ denoting the projection $\pi_i \colon M \times M \to M$ on the $i$th coordinate and $j_\varphi \colon M \to M \times M$ the graph embedding with image $\Gamma_\varphi$.
 Since $\Phi_t^\Psi$ is $C^1$-close to a symplectomorphism $\varphi$, both $\Phi_{t1}$ and $\Phi_{t2}$ are diffeomorphisms (see e.g. \cite[Exercise 3.35]{mcduffsalamon1998}). It follows that $\varphi_t= \Phi_{t2} \circ \Phi_{t1}^{-1}$ is a symplectomorphism since 
    \begin{equation}
    \omega-\varphi_t^*\omega=(\Phi_t^{\Psi})^*\big(\omega \oplus (- \omega)\big)|_{\Gamma_\varphi}=0.
    \end{equation}

We now consider a family of symplectomorphisms $\psi_t=\varphi^{-1} \circ \varphi_t$ which are $C^1$-close to the identity. Note that
\begin{equation}(\operatorname{Id} \times \varphi^{-1}) (\Gamma_{\varphi_t})=\Gamma_{\psi_t},
\end{equation}
and $\Delta_M=\Gamma_{\psi_0}$. The Hamiltonian flow $\tilde{\Phi}_t$ of $\tilde{H}_t=H_t^{\Psi} \circ (\operatorname{Id} \times \phi)$ then satisfies
\begin{equation}\tilde{\Phi}_t(\Delta_M)=\Gamma_{\psi_t}.\end{equation}
Hence, by \cite[Lemma 2.2]{akveldsalamon2001}, the derivative of the path of Lagrangian submanifolds $\Gamma_{\psi_t}$ is
\begin{equation}\label{eq: path-derivative-1}
\partial_t\Gamma_{\psi_t}=d\tilde{H}_t|_{\Gamma_{\psi_t}}.
\end{equation}
For the vector field $Y_t$ generating $\psi_t$, we compute this derivative by definition:
\begin{equation}\label{eq: path-derivative-2}
  \partial_t\Gamma_{\psi_t}=j_t^*\iota_{(0, Y_t)}\Omega=-\psi_t^*\iota_{Y_t}\omega,
\end{equation}
where $j_t \colon M \to M \times M$ is the graph embedding satisfying $j_t(x)=\big(x, \psi_t(x)\big)$. Combining~\eqref{eq: path-derivative-1} and~\eqref{eq: path-derivative-2}, we obtain
\begin{equation}\iota_{Y_t}\omega=-d(\tilde{H}_t \circ j_t \circ \psi_t^{-1}),\end{equation}
which concludes the proof.
\end{proof}

\begin{remark}
  If $\varphi$ is a finite-order Hamiltonian diffeomorphism, then it is automatically clean. As already mentioned, in \cite{OAIdeg26}, only this instance of the proposition above was used.
\end{remark}

\begin{remark}
  For the case of cleanly intersecting Lagrangian submanifolds in a given symplectic manifold, the appropriate version of the above proposition follows from the discussion in \cite[Section 3.4]{pozniak1999} (specifically \cite[Proposition 3.4.6]{pozniak1999}). We would argue that this statement is folklore among the experts.
\end{remark}

\begin{definition}
  Let $G \subset \operatorname{Ham}(M)$ be a subgroup of the group of Hamiltonian diffeomorphisms, all of whose elements are clean. We define the \emph{Arnold number} of $G$ to be
  \begin{equation}
    \operatorname{Arn}(G)\equiv\operatorname{min}_{g \in G}\mathrm{Crit}\big(\mathrm{Fix}(g)\big).
  \end{equation}
\end{definition}

Of course, Proposition~\ref{proposition: trick} implies that $\operatorname{Arn}(M) \le \operatorname{Arn}(G)$. Hence, for symplectic manifolds which have sufficiently large subgroups of clean Hamiltonian diffeomorphisms, one can try to minimize $\operatorname{Arn}(G)$ in order to get a hold of $\operatorname{Arn}(M)$.

One natural source of symplectic manifolds with large groups of Hamiltonian diffeomorphisms are flag varieties. For almost any flag variety $M=G/P$, where $G$ is a simple complex Lie group and $P$ its parabolic subgroup, it is known that the connected component $\operatorname{Aut}^0(M)$ of the identity in the group of all holomorphic automorphisms (which coincide with algebraic ones) is the group $G$ itself \cite[Chapters 3 and 4]{Akhiezer}. In most of these cases, for the natural symplectic structure on $M$, the maximal compact subgroup of $G$ acts via Hamiltonian diffeomorphisms.

\subsection{Minimizing the number of critical points for $M$ a complex flag variety}
We will constantly use the following lemma.
\begin{lemma}\label{lem: crit for flags}
  For a complex flag variety $M=G/P$, the following equality holds:
  \begin{equation}
    \mathrm{Crit}(M)=\operatorname{dim}_{\mathbb{C}}M+1.
  \end{equation}
\end{lemma}
\begin{proof}[Sketch of proof]
  The bound 
    \begin{equation}
    \mathrm{Crit}(M) \ge \operatorname{dim}_{\mathbb{C}}M+1
    \end{equation}
  is immediate since $M$ is K\"{a}hler, and 
    \begin{equation}
    \mathrm{Crit}(M) \ge \operatorname{cuplength}_{\mathbb{Q}}(M)+1.
    \end{equation}
    
  The harder part of constructing a function with the right number of critical points is based on combining the construction of a perfect Morse-Smale function of \cite{BottSamelson1958}, described in \cite[Sections 3 and 4]{DKV1983}, with \cite[Proposition 2.9]{Takens1968}. By \cite[Corollary 1.5, Corollary 3.7 and Lemma 4.2]{DKV1983}, there exists a Morse-Smale function $f \colon M \to \mathbb{R}$ with a unique minimum and a unique maximum, such that the closure of each unstable manifold for the gradient flow is an algebraic variety \cite[Lemma 4.3]{DKV1983}, implying that such a function is \emph{perfect}. Notice that all critical points have even indices for any function provided by this construction. Applying the standard result from \cite{Milnor1965}, we may replace $f$ with a function $\tilde{f}$ satisfying all the above properties, such that it is self-indexing. Now, since all critical points of $\tilde{f}$ have even index and the minimum is unique, surgery theory for Morse functions implies that each level set is connected. Hence, to each critical level set, corresponding to a critical value $0 < 2k < 2 \operatorname{dim}_{\mathbb{C}}M$, we can apply \cite[Proposition 2.9]{Takens1968} to get a function $\hat{f}$ with exactly $\operatorname{dim}_{\mathbb{C}}M+1$ critical points.
\end{proof}

\subsection{Complex Lagrangian Grassmannian manifolds}\label{section: LG}
Our first family comes from (complex) Lagrangian Grassmannians: 
    \begin{equation}
    LG(n, 2n) \cong Sp(n)/U(n) \cong Sp(2n, \mathbb{C})/P_{n},
    \end{equation}
where $P_n$ is the \emph{Siegel parabolic} subgroup consisting of elements in $Sp(2n, \mathbb{C})$ which stabilize a preferred Lagrangian subspace; it is a K\"{a}hler manifold, with K\"{a}hler structure induced from the \emph{Pl\"{u}cker embedding} $Gr(n, 2n) \hookrightarrow \mathbb{P}(\Lambda^n\mathbb{C}^{2n})$,
and $LG(n, 2n) \subset Gr(n, 2n)$ is given by the zero-set of the standard symplectic form $\omega_{\operatorname{std}}$. Moreover, this is a manifold of real dimension $\operatorname{dim} LG(n,2n) =n(n+1)$, hence   \begin{equation}\label{critLG}
\mathrm{Crit}\big(LG(n,2n)\big)=\frac{n(n+1)}{2}+1
\end{equation}
by Lemma~\ref{lem: crit for flags}.

The connected component of the identity $\operatorname{Aut}^0\big(LG(n,2n)\big)$ of the group of holomorphic automorphisms is $Sp(2n, \mathbb{C})/\{\pm1\}$ by \cite{Akhiezer}, and $Sp(n)/\{\pm 1\}$ is its maximal compact subgroup -- characterized as the group of automorphisms preserving the K\"{a}hler metric.

The action of $Sp(n)=Sp(2n, \mathbb{C}) \cap U(2n)$ is Hamiltonian, as it restricts from the Hamiltonian action of $U(2n)$ on $\mathbb{P}(\Lambda^n \mathbb{C}^{2n})$. This action is also unitary with respect to the K\"{a}hler metric on $LG(n,2n)$, as the action of $U(2n)$ is unitary on $\mathbb{P}(\Lambda^n \mathbb{C}^{2n})$ with respect to the Fubini-Study metric. Hence, $Sp(n)$ acts \emph{cleanly}, as it acts via isometries with respect to the Riemannian metric given by the real part of the K\"{a}hler metric on $LG(n,2n)$.

\begin{proposition}\label{prop: LG Arnold}
  For $Sp(n)/\{\pm 1\} \subset \operatorname{Ham}\big(LG(n,2n)\big)$, one has
  \begin{equation}\label{eq: arnold for sp}
    \operatorname{Arn}(Sp(n)/\{\pm 1\})=\begin{cases}
      \tfrac{n(n+2)}{4}+1, & n\;\mathrm{even}, \\
      \tfrac{(n+1)^2}{4}+1, & n\;\mathrm{odd}.
    \end{cases}
  \end{equation}
  The value $\operatorname{Arn}(Sp(n)/\{\pm 1\})$ is achieved by any $g \in Sp(n)$ which diagonalizes as
  \begin{equation}
    \begin{pmatrix}\label{eq: minimizer for type c}
      -\operatorname{Id}_{2\lfloor{\tfrac{n}{2}}\rfloor} & 0 \\ 0 & \operatorname{Id}_{2\lceil{\tfrac{n}{2}}\rceil}
    \end{pmatrix}.
  \end{equation}
\end{proposition}

\begin{proof}

\emph{Step 1: description of fixed points for an arbitrary $g \in Sp(n)$:} We denote by $V\cong \mathbb{C}^{2n}$ the ambient vector space. Any element $g \in Sp(n)$ is diagonalizable, and each eigenvalue $\lambda \in \mathbb{C}$ satisfies $|\lambda|=1$. Let $V=\bigoplus V_{\lambda}$ be the decomposition into $g$-eigenspaces.
  Let $\lambda_1=e^{i \theta_1}, \ldots, \lambda_k=e^{i\theta_k} \in \mathbb{C}$ be all the eigenvalues $\lambda=e^{i \theta}$ of $g$ such that $\operatorname{Im}(\lambda)>0$, or equivalently, $\theta \in (0, \pi)$. We have
  \begin{equation}
  V=\bigoplus_{i=1}^kW_{\lambda_i} \oplus V_1 \oplus V_{-1},\;\;W_\lambda=V_{\lambda}\oplus V_{\lambda^{-1}}.
  \end{equation}
  Let 
    \begin{equation}
  \operatorname{dim}_{\mathbb{C}} V_{\lambda_i}=m_i,\;\; \operatorname{dim}_{\mathbb{C}}V_{1}=2m_+,\;\;\;\mathrm{and}\;\operatorname{dim}_{\mathbb{C}}V_{-1}=2m_-,
  \end{equation}
  so $n=m_1+\cdots+m_k+m_++m_-$. We claim that
  \begin{equation}\label{eq: fixed points LG}
    \mathrm{Fix}(g)= \prod_{i=1}^k \Big(\bigsqcup_{k_i=0}^{m_i} Gr(k_i, m_i)\Big) \times LG(m_-, 2m_-) \times LG(m_+, 2m_+).
  \end{equation}
  
  Indeed, for $\lambda \neq \pm1$, the form $\omega=\omega_{\operatorname{std}}$ restricts to a nondegenerate pairing $V_\lambda \times V_{\lambda^{-1}} \to \mathbb{C}$, since $\omega(v,w)=0$ for $v \in V_{\lambda}$ and $w \in V_{\mu}$ whenever $\mu \neq \lambda^{-1}$ as $g$ preserves $\omega$. Hence, $W_\lambda\subset V$ is a symplectic subspace, and $V_{\lambda^{\pm 1}}$ are its Lagrangian subspaces (here, as before, $\lambda \neq \pm 1$). For any Lagrangian subspace $L \in \mathrm{Fix}(g)$, the restriction $L\cap W_\lambda$ is a Lagrangian subspace of $W_\lambda$. We notice that the subspace $L_\lambda= L \cap V_\lambda$ completely determines $L_{\lambda^{-1}}=L \cap V_{\lambda^{-1}}$ as 
  \begin{equation}
  L_{\lambda^{-1}}=\big\{v \in V_{\lambda^{-1}}: \omega(v, L_\lambda)=0\big\}.
  \end{equation}
  Thus, fixing the dimension $0 \le k_i \le m_i$ of the restriction $L_{\lambda_i}$ gives rise to the component $Gr(k_i, m_i)$ in a potential product as in~\eqref{eq: fixed points LG}.
  For the remaining two potential eigenvalues, the whole subspace $V_{\pm 1}$ is symplectic, and therefore $L_{\pm 1}=L \cap V_{\pm 1}$ is a Lagrangian subspace in it, which gives rise to the $LG(m_{\pm}, 2m_{\pm})$.

\emph{Step 2: ``stretching films'' between fixed loci:}
  Applying Lemma~\ref{lem: crit for flags}, we see that, for each connected component $Z$ of $\mathrm{Fix}(g)$, we have $\mathrm{Crit}(Z)=\operatorname{dim}_{\mathbb{C}}Z+1$. Intuitively, the sum of dimensions over the connected components is smaller for those $g$ which have fewer connected components. We implement this idea through a ``film stretching'' procedure.
  
  Let $g_j^{\pm}$, $1 \le j\le k$, be an element of $Sp(n)$ obtained by ``rotating" the full block \begin{equation}\label{eq: block-rotation}
  \begin{pmatrix} e^{i \theta_j}\operatorname{Id}_{m_j} & 0 \\ 0 & e^{-i \theta_j} \operatorname{Id}_{m_j}
  \end{pmatrix}\end{equation}
  of $g$, corresponding to its action on $W_{\lambda_j}$, into the block \begin{equation}\begin{pmatrix} \pm\operatorname{Id}_{m_j} & 0 \\ 0 & \pm \operatorname{Id}_{m_j}
  \end{pmatrix}.
  \end{equation}
The effect of this rotation is the replacement of the term \begin{equation}
    Z=\Big(\bigsqcup_{k_j=0}^{m_j} Gr(k_j, m_j)\Big) \times LG(m_{\pm}, 2m_{\pm})
  \end{equation} in~\eqref{eq: fixed points LG} with 
    \begin{equation}
    Z_j=LG\big(m_j+m_{\pm}, 2(m_j+m_{\pm})\big).
    \end{equation}
    One can think of $Z_j$ as a film stretched between components of $Z$.
  
  Now, we simply count sums of the dimensions shifted by $1$ over the connected components constituting $Z$ and $Z_j$:
  \begin{align}
    \mathrm{Crit}(Z)&= \sum_{k_j=0}^{m_j}\Big[\binom{m_{\pm}+1}{2} + k_j(m_j-k_j)+1 \Big] \nonumber \\
    &=( m_j+1) \binom{m_{\pm}+1}{2}+ \binom{m_j+1}{3}+(m_j+1), \\
    \mathrm{Crit} Z_j&=\binom{m_j+m_{\pm}+1}{2}+1 \nonumber \\
    &=\binom{m_j+1}{2}+\binom{m_\pm+1}{2}+m_jm_\pm+1.
  \end{align}
Hence, we have
  \begin{align}
    \mathrm{Crit}(Z)-\mathrm{Crit}(Z_j)&=\Big[ \binom{m_j+1}{3}+m_j-\binom{m_j+1}{2}\Big]+\Big[ m_j\binom{m_{\pm}+1}{2}-m_jm_{\pm}\Big] \nonumber \\
    &=\binom{m_j}{3}+m_j \binom{m_{\pm}}{2},
  \end{align}
  which is nonnegative whenever $m_j, m_{\pm} \ge 0$. From this computation, it follows immediately that 
  \begin{equation*}
    \mathrm{Crit}\big(\mathrm{Fix}(g_j^{\pm})\big) \le \mathrm{Crit}\big(\mathrm{Fix}(g)\big).
  \end{equation*}
  
\begin{remark}
Let us note that this inequality can be obtained more conceptually by observing that there is a Hamiltonian action of a $1$-dimensional torus by elements of the shape~\eqref{eq: block-rotation} on $Z_j$ for which $Z$ is exactly the fixed point locus. This argument is outlined in the proof of Proposition~\ref{prop: arn for type B}.
\end{remark}

\emph{Step 3: optimization:}
  By the previous step, we may assume that $k=0$ and $g$ only has eigenvalues $\pm1$, i.e., $n=m_++m_-$, and
  \begin{equation}
    \mathrm{Fix}(g)=LG(m_-, 2m_-) \times LG(m_+, 2m_+)
  \end{equation}
  with \begin{equation}
    \mathrm{Crit}\big(\mathrm{Fix}(g)\big)=\binom{m_-+1}{2}+\binom{m_++1}{2}+1;
  \end{equation}
    as a function of $m_- \in \mathbb{Z}$, with $0 \le m_- \le n$, this takes its minimal value at $m_-=\lfloor{\tfrac{n}{2}} \rfloor$, and is given by~\eqref{eq: arnold for sp}.
\end{proof}
Combining Propositions~\ref{proposition: trick} and \ref{prop: LG Arnold}, we immediately obtain the following.
\begin{corollary}\label{cor: arn for LG}
  For $M=LG(n,2n)$, the Arnold number $\operatorname{Arn}(M)$ is bounded above by ~\eqref{eq: arnold for sp}, which is approximately two times less than $\mathrm{Crit}(M)$ given in \eqref{critLG}.
\end{corollary}
\begin{remark}
  Lagrangian Grassmannians belong to the family of type $C$ isotropic Grassmannians $IG(k, 2n)$ (and more generally type $C$ flag varieties). Determining the Arnold number for these more general manifolds is also of great interest. Our findings, which we do not report here, suggest that $\operatorname{Arn}(Sp(n)/\{\pm 1\})=\mathrm{Crit}\big(IG(k, 2n)\big)$ for $k \neq 1, n$. The main reason for this is that, for a $k$-isotropic subspace $L$ with $k <n$, the dimension of the intersection of $L$ with an arbitrary symplectic subspace is not predetermined (like in the case of Lagrangians); this causes non-identity elements in $Sp(n)$ to have numerous connected components.
\end{remark}

\begin{remark}
  One may hope that establishing the sharp upper bound for the Arnold number in the case of an arbitrary symplectic manifold $M$ may go along the lines of ``stretching films'' between fixed components of a given Hamiltonian diffeomorphism at each step of induction.
\end{remark}

\subsection{Isotropic Grassmannian manifolds of type $B$}
The discussion here is parallel to that in Section~\ref{section: LG}, so we keep it more confined. Let $B$ be a nondegenerate symmetric bilinear form on $V\cong \mathbb{C}^{2n+1}$. We consider
\begin{equation}
  OG(n-k, 2n+1) =\big\{W \in Gr(n -k, 2n+1) : B|_W \equiv 0\big\};
\end{equation}
it is a K\"{a}hler manifold with a group of holomorphic automorphisms preserving the K\"{a}hler metric, and the connected component of the identity is $SO(2n+1)$ if $k>0$; for $k=0$, it contains this group. Each element $g \in SO(2n+1)$ is diagonalizable and acts cleanly on $OG(n-k, 2n+1)$. We note that $\mathrm{Crit}\big(OG(n-k, 2n+1)\big)$ is given in Table \ref{tab:arn-bounds}.

\begin{proposition}\label{prop: arn for type B}
  For $M=OG(n-k, 2n+1)$, $k>0$, and $G=SO(2n+1) \subset \operatorname{Ham}(M)$, one has
  \begin{equation}\label{eq: arn for type B}
    \operatorname{Arn}(G)=\mathrm{Crit}(M)-(n-k)=(n-k)(n+k+1)-\binom{n-k+1}{2}-(n-k)+1.
  \end{equation}
  In the exceptional case $k=0$, one has
  \begin{equation}\label{eq: optimal for k=0 type b}
    \operatorname{Arn}(G)= \lfloor\tfrac{n^2}{2} \rfloor +2.
  \end{equation}
\end{proposition}
\begin{proof}[Sketch of proof]
\emph{Step 1: description of the element minimizing the Arnold number:}
  For $k>0$, the minimum is attained on 
  \begin{equation}g= \begin{pmatrix}
    -\operatorname{Id}_{2n} & 0 \\ 0 & 1
  \end{pmatrix}\end{equation}
  since $\mathrm{Fix}(g)=OG(n-k, 2n)$. Indeed, $g$ acts on an isotropic subspace $V$ by reflecting it in the plane orthogonal to the last basis vector $e_{2n+1}$, and since $B(e_{2n+1}, e_{2n+1})=1$, the subspace $V$ is fixed under such a reflection only if it is contained in the orthogonal complement $\langle e_{2n+1} \rangle^{\perp}$.

For $k=0$ one instead takes \begin{equation*}g= \begin{pmatrix}
    -\operatorname{Id}_{2\lceil\frac{n+1}{2}\rceil} & 0 \\ 0 & \operatorname{Id}_{2\lfloor\frac{n-1}{2}\rfloor+1}
  \end{pmatrix}\end{equation*}
  with \begin{equation}
  \operatorname{Fix}(g)=OG\Big(\lceil\frac{n+1}{2}\rceil, 2\lceil\frac{n+1}{2}\rceil\Big) \times OG\Big(\lceil\frac{n+1}{2}\rceil, 2\lceil\frac{n+1}{2}\rceil+1\Big).
  \end{equation}

\emph{Step 2: film stretching.}
  The stretching argument is similar to the one performed in the proof of Proposition~\ref{prop: LG Arnold}; here, we proceed as follows. Up to conjugation, an element $g \in SO(2n+1)$ has the form
  \begin{equation}
  g= \operatorname{diag}(e^{i\theta_1}\operatorname{Id}_{m_1}, e^{-i\theta_1}\operatorname{Id}_{m_1}, \dots, e^{i\theta_\ell}\operatorname{Id}_{m_\ell}, e^{-i\theta_\ell}\operatorname{Id}_{m_\ell}, -\operatorname{Id}_{2m_-}, \operatorname{Id}_{2m_++1})
  \end{equation}
  with $m_1+\dots+m_\ell+m_-+m_+=n$ and $\theta_1, \ldots, \theta_{\ell} \in (0, \pi)$. Let 
  \begin{equation}
  g_1^-=\operatorname{diag}(\operatorname{-Id}_{2m_1}, e^{i\theta_2}\operatorname{Id}_{m_2}, e^{-i\theta_2}\operatorname{Id}_{m_2}, \dots,e^{i\theta_\ell}\operatorname{Id}_{m_\ell}, e^{-i\theta_\ell}\operatorname{Id}_{m_\ell}, -\operatorname{Id}_{2m_-}, \operatorname{Id}_{2m_++1}).
  \end{equation}
  
  Let $W=V_{e^{i\theta_1}}^g \oplus V_{e^{-i\theta_1}}^g$, where $V_{e^{\pm i\theta_1}}^g$ is the eigenspace of $g$ for the eigenvalue $\lambda=e^{\pm i \theta_1}$. For each $ 0\le s \le m_1+m_-$, we compare the components of $\mathrm{Fix}(g)$ and $\mathrm{Fix}(g_1^-)$ consisting of those isotropic subspaces $U$ such that $\operatorname{dim}U \cap V_{-1}^{g_1^-}=s$, where $V_{-1}^{g_1^-}=W \oplus V_{-1}$ is the eigenspace of $g_1^{-}$ for the eigenvalue $-1$. The restriction $U \cap {V_{-1}^{g_1^-}}$ of $U \in \mathrm{Fix}(g_1^-)$ to $V_{-1}^{g_1^-}$ can be any isotropic subspace: \begin{equation}\label{eq: component for rotated g}
    U \cap {V_{-1}^{g_1^-}} \in OG\big(s, 2(m_1+m_-)\big).
  \end{equation} 
  For $U$ to belong to $\mathrm{Fix}(g)$, the conditions are stricter, as both restrictions 
  \begin{equation}
  U \cap V^g_{e^{i \theta_1}} \text{ and }U \cap V^g_{e^{-i \theta_1}}
  \end{equation}
  have to be isotropic and orthogonal to each other, i.e., 
  \begin{gather}\label{eq: components for closer g}
  U \cap V_{-1}^{g_1^-} \in \bigsqcup_{a+b+c=s, \, a+b \le m_1} Z_{a, b} \times OG(c, 2m_-),\;\mathrm{where} \\
  Z_{a,b}=\big\{(U_1, U_2) \in Gr(a, V^g_{e^{i \theta_1}}) \times Gr(b, V_{e^{-i\theta_1}}^g) : B(U_1, U_2)=0\big\}. \nonumber
  \end{gather}
    This argument shows that $\mathrm{Fix}(g) \subset \mathrm{Fix}(g_1^-)$ by comparing~\eqref{eq: component for rotated g} and~\eqref{eq: components for closer g}, and one can treat $\mathrm{Fix}(g_1^-)$ as a film stretched between components of $\mathrm{Fix}(g)$.
 
  To compare the minimal numbers of critical points on these two sets, we consider the Hamiltonian action of a $1$-dimensional torus $T \cong S^1$ with $\theta \in T$ acting via multiplication by
  \begin{equation}g_\theta=(e^{i\theta} \operatorname{Id}_{m_1}, e^{-i\theta} \operatorname{Id}_{m_1}, 1, \dots, 1) \in SO(2n+1).\end{equation}
  Notice that, since $g_\theta$ commutes with $g_1^-$, it induces a Hamiltonian action on $\mathrm{Fix}(g_1^-)$, while $\mathrm{Fix}(g)$ coincides with the fixed locus of the whole torus $T$ on $\mathrm{Fix}(g_1^-)$. Hence, $\mathrm{Fix}(g)$ coincides with the critical locus of the moment map $\mu \colon \mathrm{Fix}(g_1^-) \to \mathbb{R}$, and the moment map has Morse-Bott critical loci. 
  
  Applying an argument similar to the one in \cite[Section 5]{BanyagaHurtubise2009}, one can show \begin{equation}
  \mathrm{Crit}\big(\mathrm{Fix}(g^{-}_1)\big) \le \mathrm{Crit} \big(\mathrm{Fix}(g)\big),
  \end{equation}
  since $\mathrm{Fix}(g)$ is a Morse-Bott critical locus of a Morse-Bott function $\mu \colon \mathrm{Fix}(g_1^{-}) \to \mathbb{R}$. To do so, one simply perturbs $\mu$ by $\varepsilon \chi \cdot \tilde{f}$, where $\tilde{f}$ is an extension to a tubular neighborhood of $\mathrm{Fix}(g) \subset \mathrm{Fix}(g_1^-)$ of a smooth function $f \in C^\infty (\mathrm{Fix}(g))$ with \begin{equation}
  \mathrm{Crit}(f)=\mathrm{Crit}\big(\mathrm{Fix}(g)\big),
  \end{equation}and $\chi$ is an appropriate bump function. Using the Morse-Bott lemma, it is immediate to see that, for small $\varepsilon>0$, the function $\mu+\varepsilon \chi\cdot\tilde{f}$ has exactly $\mathrm{Crit}(f)$ critical points.

  \emph{Step 3: optimization:}
  In particular, we may assume that 
    \begin{equation}
    g= \begin{pmatrix}
    -\operatorname{Id}_{2m_-} & 0 \\ 0 & \operatorname{Id}_{2m_++1}
  \end{pmatrix}\end{equation}
  with $m_-+m_+=n$. We then have 
  \begin{equation}
    \mathrm{Fix}(g)=\bigsqcup_{j=\max(0, m_- -k)}^{\min{(m_-, n-k)}}OG(j, 2m_-) \times OG(n-k-j, 2m_++1),
  \end{equation}
  and, by Proposition~\ref{proposition: trick}, it is straightforward to check that $\mathrm{Crit}\big(\mathrm{Fix}(g)\big)$ is minimal for $m_+=0$ whenever $k>0$, for which it is given by~\eqref{eq: arn for type B}. 
  
  The case $k=0$ is special since $OG(n, 2n)$ is disconnected, and $m_+=0$ is not optimal. In this case, we have
  \begin{equation}\mathrm{Fix}(g)= OG(m_-, 2m_-) \times OG(m_+, 2m_++1),\end{equation}
  and $\mathrm{Crit}\big(\mathrm{Fix}(g)\big)$ is minimal for $m_-= \lceil{\tfrac{n+1}{2}}\rceil$, for which it is given by~\eqref{eq: optimal for k=0 type b}.
\end{proof}

Combining Propositions~\ref{proposition: trick} and \ref{prop: arn for type B}, we immediately obtain the following

\begin{corollary}\label{cor: arn for type B}
  For $M=OG(n-k, 2n+1)$, $k>0$, and $G=SO(2n+1) \subset \operatorname{Ham}(M)$, one has
  \begin{equation}
    \operatorname{Arn}(M) \le\mathrm{Crit}(M)-(n-k).
  \end{equation}
  In the exceptional case $k=0$, one has
  \begin{equation}
    \operatorname{Arn}(M) \le \lfloor\tfrac{n^2}{2} \rfloor +2.
  \end{equation}
\end{corollary}

\begin{remark}
  We briefly note that in type $D$, i.e., for Grassmannians $OG(n-k, 2n)$, the situation is similar to that of type $C$. For the maximal isotropic Grassmannian $OG^+(n,2n)$ (i.e., for the spinor variety), the matrix~\eqref{eq: minimizer for type c} also provides improvements as compared with $\mathrm{Crit}\big(OG(n, 2n)\big)$, although not as striking as in type $C$ since there are two connected components in the fixed locus for its action on $OG(n, 2n)$. The case $OG(1, 2n)$ is discussed in the sequel.
\end{remark}

\subsection{Quadrics}\label{section: quadrics}
Let us note that the $n$-dimensional quadric $Q^n$ is precisely $OG(1,n+2)$, and $\operatorname{Arn}(Q^n) \ge n$ is well-known by \cite[Corollary 1.3]{Schwarz1998QuantumCupLength}. For $n=2$, the example of a Hamiltonian diffeomorphism on $Q^2\cong \mathbb{P}^1 \times \mathbb{P}^1$ with exactly $n=2$ fixed points is provided in \cite{jiao2026twofixed}, and the example for odd $n$ is given in the previous section; the latter upgrades the construction from \cite{OAIdeg26} for the case $n=3$. Here, we explain how to get a Hamiltonian diffeomorphism on any quadric $Q^n$ with $n\ge 6$ and even.

\begin{proposition}\label{prop: quadrics}
  For any quadric $Q^n$ with $n \ge 2$, $n \neq 4$, one has
  \begin{equation}
    \operatorname{Arn}(Q^n)=n.
  \end{equation}
\end{proposition}

\begin{proof}
  Since the case of $n$ odd is already proved in Corollary~\ref{cor: arn for type B}, we concentrate our attention on the case of even $n$ with $n\ge 6$. We consider the action of \begin{equation}
  g =\begin{pmatrix}
    -\operatorname{Id}_4 & 0 \\
    0 & \operatorname{Id}_{n-2}
  \end{pmatrix}\end{equation}
  whose fixed point set $\mathrm{Fix}(g)$ is given by
  \begin{equation}
    OG(1, 4) \sqcup OG(1, n-2) \cong Q^2 \sqcup Q^{n-4},
  \end{equation}
  with both components connected by the assumption on $n$. Then, since
  \begin{equation}
  \mathrm{Crit}\big(\mathrm{Fix}(g)\big)=\mathrm{Crit}(Q^2)+\mathrm{Crit}(Q^{n-4})=n,
  \end{equation}
  the claim now follows by applying Proposition~\ref{proposition: trick}.
\end{proof}

\begin{remark}
  The reason our method does not work in the case $n=4$ with $Q^n \cong Gr(2,4)$ is due to the fact that $OG(1,n-2) \cong \operatorname{pt} \sqcup \operatorname{pt}$ is disconnected. We expect that one can build a Hamiltonian diffeomorphism on $Gr(2,4)$ with exactly $4$ fixed points, but we do not pursue it here.
\end{remark}
\bibliography{References}{}

@article{Flo88b,
  author = {Floer, A.},
  title = {{Morse theory for Lagrangian intersections}},
  journal = {Journal of Differential Geometry},
  volume = {28},
  number = {3},
  pages = {513--547},
  year = {1988},
}

@article{Flo89a,
  author = {Floer, A.},
  title = {{Cuplength estimates on Lagrangian intersections}},
  journal = {Communications on Pure and Applied Mathematics},
  volume = {42},
  number = {4},
  pages = {335--356},
  year = {1989},
}

@article{Flo89b,
  author = {Floer, A.},
  title = {{Symplectic fixed points and holomorphic spheres}},
  journal = {Communications in Mathematical Physics},
  volume = {120},
  number = {4},
  pages = {575--611},
  year = {1989},
}

@article{Hof85,
  author = {Hofer, H.},
  title = {{Lagrangian embeddings and critical point theory}},
  journal = {Annales de l'Institut Henri Poincar\'e. Analyse non lin\'eaire},
  volume = {2},
  number = {6},
  pages = {407--462},
  year = {1985},
}

@article{Hof88,
  author = {Hofer, H.},
  title = {{Lusternik-Schnirelman-theory for Lagrangian intersections}},
  journal = {Annales de l'Institut Henri Poincar\'e. Analyse non lin\'eaire},
  volume = {5},
  number = {5},
  pages = {465--499},
  year = {1988},
}

@unpublished{AB21,
  author = {Abouzaid, M. and Blumberg, A.},
  title = {{Arnold conjecture and Morava K-theory}},
  year = {2021},
  note = {arXiv:2103.01507},
}

@unpublished{BX22,
  author = {Bai, S. and Xu, G.},
  title = {{Arnold conjecture over integers}},
  year = {2022},
  note = {arXiv:2209.08599},
}

@unpublished{Rez22,
  author = {Rezchikov, S.},
  title = {{Integral Arnol'd conjecture}},
  year = {2022},
  note = {arXiv:2209.11165},
}

@article{akveldsalamon2001,
  author = {Akveld, M. and Salamon, D.},
  title = {{Loops of Lagrangian submanifolds and pseudoholomorphic discs}},
  journal = {Geometric and Functional Analysis},
  volume = {11},
  number = {4},
  pages = {609--650},
  year = {2001},
  doi = {10.1007/PL00001680},
  url = {https://doi.org/10.1007/PL00001680},
  issn = {1016-443X,1420-8970},
  fjournal = {Geometric and Functional Analysis},
  mrclass = {53D45 (37J05 53D12)},
  mrnumber = {1866797},
  mrreviewer = {Francisco\ Presas},
}

@book{mcduffsalamon1998,
  author = {McDuff, D. and Salamon, D.},
  title = {{Introduction to symplectic topology}},
  series = {Oxford Mathematical Monographs},
  edition = {Second},
  pages = {x+486},
  publisher = {The Clarendon Press, Oxford University Press, New York},
  year = {1998},
  isbn = {0-19-850451-9},
  mrclass = {53D35 (53D40 57R17 57R57 57R58)},
  mrnumber = {1698616},
  mrreviewer = {Hansj\"org\ Geiges},
}

@incollection{pozniak1999,
  author = {Po{\'z}niak, M.},
  title = {{Floer homology, Novikov rings and clean intersections}},
  booktitle = {Northern {C}alifornia {S}ymplectic {G}eometry {S}eminar},
  series = {Amer. Math. Soc. Transl. Ser. 2},
  volume = {196},
  pages = {119--181},
  publisher = {Amer. Math. Soc., Providence, RI},
  year = {1999},
  doi = {10.1090/trans2/196/08},
  url = {https://doi.org/10.1090/trans2/196/08},
  isbn = {0-8218-2075-3},
  mrclass = {53D40 (57R58)},
  mrnumber = {1736217},
  mrreviewer = {David\ E.\ Hurtubise},
}

@article{Arnold1965surune,
  author = {Arnold, V.},
  title = {{Sur une propri\'et\'e{} topologique des applications globalement canoniques de la m\'ecanique classique}},
  journal = {Comptes Rendus Hebdomadaires des S\'eances de l'Acad\'emie des Sciences},
  volume = {261},
  pages = {3719--3722},
  year = {1965},
  issn = {0001-4036},
  fjournal = {Comptes Rendus Hebdomadaires des S\'eances de l'Acad\'emie des Sciences},
  mrclass = {57.50 (34.65)},
  mrnumber = {193645},
  mrreviewer = {E.\ A.\ Feldman},
}

@book{Arnold1978methods,
  author = {Arnold, V.},
  title = {{Mathematical methods of classical mechanics}},
  series = {Graduate Texts in Mathematics},
  volume = {60},
  pages = {x+462},
  publisher = {Springer-Verlag, New York-Heidelberg},
  year = {1978},
  note = {Translated from the Russian by K. Vogtmann and A. Weinstein},
  isbn = {0-387-90314-3},
  mrclass = {58F05 (70.58)},
  mrnumber = {690288},
  mrreviewer = {J.\ S.\ Joel},
}

@unpublished{Eliashberg1978,
  author = {Eliashberg, Y.},
  title = {{An estimate of the number of fixed points of transformations preserving area}},
  year = {1978},
  note = {preprint, Syktyvkar},
}

@misc{Weinstein1983,
  author = {Weinstein, A.},
  title = {{$C^0$ perturbation theorems for symplectic fixed points and Lagrangian intersections}},
  year = {1984},
  note = {South Rhone Seminar on Geometry, III (Lyon, 1983). Hermann, Paris: Travaux en Cours},
}

@article{Fortune1985,
  author = {Fortune, B.},
  title = {{A symplectic fixed point theorem for $\mathbf{C}\mathrm{P}^n$}},
  journal = {Inventiones Mathematicae},
  volume = {81},
  number = {1},
  pages = {29--46},
  year = {1985},
  doi = {10.1007/BF01388770},
  url = {https://doi.org/10.1007/BF01388770},
  issn = {0020-9910,1432-1297},
  fjournal = {Inventiones Mathematicae},
  mrclass = {58F05 (49A40 58C30)},
  mrnumber = {796189},
  mrreviewer = {Frans\ Cantrijn},
}

@misc{jiao2026twofixed,
  author = {Jiao, H.},
  title = {{A smooth Hamiltonian diffeomorphism with two fixed points on $S^2 \times S^2$}},
  howpublished = {arXiv:2609.33626},
  year = {2026},
  primaryclass = {math.SG},
}

@article{Floer1986,
  author = {Floer, A.},
  title = {{Proof of the Arnol'd conjecture for surfaces and generalizations to certain K\"ahler manifolds}},
  journal = {Duke Mathematical Journal},
  volume = {53},
  number = {1},
  pages = {1--32},
  year = {1986},
  doi = {10.1215/S0012-7094-86-05301-9},
  url = {https://doi.org/10.1215/S0012-7094-86-05301-9},
  issn = {0012-7094,1547-7398},
  fjournal = {Duke Mathematical Journal},
  mrclass = {58F22 (53C55 53C57 58F05)},
  mrnumber = {835793},
  mrreviewer = {A.\ Morimoto},
}

@article{CZ1983,
  author = {Conley, C. C. and Zehnder, E.},
  title = {{The Birkhoff-Lewis fixed point theorem and a conjecture of V. I. Arnol'd}},
  journal = {Inventiones Mathematicae},
  volume = {73},
  number = {1},
  pages = {33--49},
  year = {1983},
  doi = {10.1007/BF01393824},
  url = {https://doi.org/10.1007/BF01393824},
  issn = {0020-9910,1432-1297},
  fjournal = {Inventiones Mathematicae},
  mrclass = {58F05 (58C30)},
  mrnumber = {707347},
  mrreviewer = {Yu.\ E.\ Gliklikh},
}

@article{Oh1990symplectic,
  author = {Oh, Y.-G.},
  title = {{A symplectic fixed point theorem on $T^{2n} \times \mathbf{C}\mathrm{P}^k$}},
  journal = {Mathematische Zeitschrift},
  volume = {203},
  number = {4},
  pages = {535--552},
  year = {1990},
  doi = {10.1007/BF02570755},
  url = {https://doi.org/10.1007/BF02570755},
  issn = {0025-5874,1432-1823},
  fjournal = {Mathematische Zeitschrift},
  mrclass = {58F05 (58C30 58F20)},
  mrnumber = {1044063},
  mrreviewer = {J.\ W.\ Robbin},
}

@incollection{LeOno1994,
  author = {L\^e, H. V. and Ono, K.},
  title = {{Cup-length estimates for symplectic fixed points}},
  booktitle = {Contact and symplectic geometry ({C}ambridge, 1994)},
  series = {Publ. Newton Inst.},
  volume = {8},
  pages = {268--295},
  publisher = {Cambridge Univ. Press, Cambridge},
  year = {1996},
  isbn = {0-521-57086-7},
  mrclass = {58F05 (58E05 58F20)},
  mrnumber = {1432466},
  mrreviewer = {Karl\ Friedrich\ Siburg},
}

@article{LS1985,
  author = {Laudenbach, F. and Sikorav, J.-C.},
  title = {{Persistance d'intersection avec la section nulle au cours d'une isotopie hamiltonienne dans un fibr\'e{} cotangent}},
  journal = {Inventiones Mathematicae},
  volume = {82},
  number = {2},
  pages = {349--357},
  year = {1985},
  doi = {10.1007/BF01388807},
  url = {https://doi.org/10.1007/BF01388807},
  issn = {0020-9910,1432-1297},
  fjournal = {Inventiones Mathematicae},
  mrclass = {58F05},
  mrnumber = {809719},
  mrreviewer = {I.\ Vaisman},
}

@article{AlbersHein2016,
  author = {Albers, P. and Hein, D.},
  title = {{Cuplength estimates in Morse cohomology}},
  journal = {Journal of Topology and Analysis},
  volume = {8},
  number = {2},
  pages = {243--272},
  year = {2016},
  doi = {10.1142/S1793525316500102},
  url = {https://doi.org/10.1142/S1793525316500102},
  issn = {1793-5253,1793-7167},
  fjournal = {Journal of Topology and Analysis},
  mrclass = {57R58 (55N45 57R70)},
  mrnumber = {3474075},
  mrreviewer = {Darko\ Milinkovi\'c},
}

@article{AsselleStarostka26,
  author = {Asselle, L. and Starostka, M.},
  title = {{On the degenerate Arnold conjecture on $\mathbb{T}^{2M} \times \mathbb{CP}^{N}$}},
  journal = {Topological Methods in Nonlinear Analysis},
  volume = {67},
  number = {1},
  pages = {19--34},
  year = {2026},
  issn = {1230-3429},
  fjournal = {Topological Methods in Nonlinear Analysis},
  mrclass = {37J11 (37B30 58E05)},
  mrnumber = {5074720},
}

@incollection{chaperon1983,
  author = {Chaperon, M.},
  title = {{Quelques questions de g\'eom\'etrie symplectique (d'apr\`es, entre autres, Poincar\'e, Arnol'd, Conley et Zehnder)}},
  booktitle = {Bourbaki seminar, {V}ol. 1982/83},
  series = {Ast\'erisque},
  volume = {105--106},
  pages = {231--249},
  publisher = {Soc. Math. France, Paris},
  year = {1983},
  mrclass = {58F05 (58-02)},
  mrnumber = {728991},
  mrreviewer = {Richard\ C.\ Swanson},
}

@article{AAC2025,
  author = {Alizadeh, H. and Atallah, M. S. and Cant, D.},
  title = {{Lagrangian intersections and the spectral norm in convex-at-infinity symplectic manifolds}},
  journal = {Mathematische Zeitschrift},
  volume = {310},
  number = {1},
  pages = {Article 2, 48 pp.},
  year = {2025},
  doi = {10.1007/s00209-025-03710-0},
  url = {https://doi.org/10.1007/s00209-025-03710-0},
  issn = {0025-5874,1432-1823},
  fjournal = {Mathematische Zeitschrift},
  mrclass = {53D12 (37J06 53D40)},
  mrnumber = {4879080},
  mrreviewer = {Weiyi\ Zhang},
}

@article{BHS2018,
  author = {Buhovsky, L. and Humili\`ere, V. and Seyfaddini, S.},
  title = {{A {$C^0$} counterexample to the Arnold conjecture}},
  journal = {Inventiones Mathematicae},
  volume = {213},
  number = {2},
  pages = {759--809},
  year = {2018},
  doi = {10.1007/s00222-018-0797-x},
  url = {https://doi.org/10.1007/s00222-018-0797-x},
  issn = {0020-9910,1432-1297},
  fjournal = {Inventiones Mathematicae},
  mrclass = {53D35},
  mrnumber = {3827210},
  mrreviewer = {Chun-Gen\ Liu},
}

@incollection{audin2014,
  author = {Audin, M.},
  title = {{Vladimir Igorevich Arnold and the invention of symplectic topology}},
  booktitle = {Contact and symplectic topology},
  series = {Bolyai Soc. Math. Stud.},
  volume = {26},
  pages = {1--25},
  publisher = {J\'anos Bolyai Math. Soc., Budapest},
  year = {2014},
  doi = {10.1007/978-3-319-02036-5_1},
  url = {https://doi.org/10.1007/978-3-319-02036-5_1},
  isbn = {978-3-319-02035-8; 978-3-319-02036-5},
  mrclass = {01A70 (53-03)},
  mrnumber = {3220939},
}

@article{FO99,
  author = {Fukaya, K. and Ono, K.},
  title = {{Arnold conjecture and Gromov-Witten invariant}},
  journal = {Topology},
  volume = {38},
  number = {5},
  pages = {933--1048},
  year = {1999},
  doi = {10.1016/S0040-9383(98)00042-1},
  url = {https://doi.org/10.1016/S0040-9383(98)00042-1},
  issn = {0040-9383},
  fjournal = {Topology. An International Journal of Mathematics},
  mrclass = {53D45 (37J10 37J45 53D40 57R17 57R58)},
  mrnumber = {1688434},
  mrreviewer = {David\ E.\ Hurtubise},
}

@article{FW22,
  author = {Filippenko, B. and Wehrheim, K.},
  title = {{A polyfold proof of the Arnold conjecture}},
  journal = {Selecta Mathematica},
  volume = {28},
  number = {1},
  pages = {Article 11, 73 pp.},
  year = {2022},
  doi = {10.1007/s00029-021-00680-z},
  url = {https://doi.org/10.1007/s00029-021-00680-z},
  issn = {1022-1824,1420-9020},
  fjournal = {Selecta Mathematica. New Series},
  mrclass = {53D05 (37J06 46T20)},
  mrnumber = {4346508},
  mrreviewer = {Guangbo\ Xu},
}

@article{LT98,
  author = {Liu, G. and Tian, G.},
  title = {{Floer homology and Arnold conjecture}},
  journal = {Journal of Differential Geometry},
  volume = {49},
  number = {1},
  pages = {1--74},
  year = {1998},
  url = {http://projecteuclid.org/euclid.jdg/1214460936},
  issn = {0022-040X,1945-743X},
  fjournal = {Journal of Differential Geometry},
  mrclass = {58E05 (55N35 57R15 57R70 58F05)},
  mrnumber = {1642105},
  mrreviewer = {Jean-Claude\ Sikorav},
}

@article{Pardon16,
  author = {Pardon, J.},
  title = {{An algebraic approach to virtual fundamental cycles on moduli spaces of pseudo-holomorphic curves}},
  journal = {Geometry \& Topology},
  volume = {20},
  number = {2},
  pages = {779--1034},
  year = {2016},
  doi = {10.2140/gt.2016.20.779},
  url = {https://doi.org/10.2140/gt.2016.20.779},
  issn = {1465-3060,1364-0380},
  fjournal = {Geometry \& Topology},
  mrclass = {53D35 (37J10 53D37 53D40 53D42 53D45 54B40 57R17)},
  mrnumber = {3493097},
  mrreviewer = {Sonja\ Hohloch},
}

@misc{OAIdeg26,
  author = {{OpenAI}},
  title = {{Three fixed points on the symplectic quadric threefold}},
  howpublished = {OpenAI Math Release preprint},
  year = {2026},
  url = {https://github.com/openai/math/blob/main/preprints/A-degenerate-counterexample-to-the-critical-number-Arnold-bound-September-23-2026/paper.pdf},
}

@article{RO1999,
  author = {Rudyak, Y. B. and Oprea, J.},
  title = {{On the Lusternik-Schnirelmann category of symplectic manifolds and the Arnold conjecture}},
  journal = {Mathematische Zeitschrift},
  volume = {230},
  number = {4},
  pages = {673--678},
  year = {1999},
  doi = {10.1007/PL00004709},
  url = {https://doi.org/10.1007/PL00004709},
  issn = {0025-5874,1432-1823},
  fjournal = {Mathematische Zeitschrift},
  mrclass = {53D40 (55M30)},
  mrnumber = {1686579},
  mrreviewer = {David\ E.\ Hurtubise},
}

@article{BottSamelson1958,
  author = {Bott, R. and Samelson, H.},
  title = {{Applications of the theory of Morse to symmetric spaces}},
  journal = {American Journal of Mathematics},
  volume = {80},
  pages = {964--1029},
  year = {1958},
  doi = {10.2307/2372843},
  url = {https://doi.org/10.2307/2372843},
  issn = {0002-9327,1080-6377},
  fjournal = {American Journal of Mathematics},
  mrclass = {55.00 (22.00)},
  mrnumber = {105694},
  mrreviewer = {C.\ T.\ Yang},
}

@book{Akhiezer,
  author = {Akhiezer, D. N.},
  title = {{Lie group actions in complex analysis}},
  series = {Aspects of Mathematics},
  volume = {E27},
  pages = {viii+201},
  publisher = {Friedr. Vieweg \& Sohn, Braunschweig},
  year = {1995},
  doi = {10.1007/978-3-322-80267-5},
  url = {https://doi.org/10.1007/978-3-322-80267-5},
  isbn = {3-528-06420-X},
  mrclass = {32M05 (22-02 22E10 32M10)},
  mrnumber = {1334091},
  mrreviewer = {B.\ Gilligan},
}

@article{DKV1983,
  author = {Duistermaat, J. J. and Kolk, J. A. C. and Varadarajan, V. S.},
  title = {{Functions, flows and oscillatory integrals on flag manifolds and conjugacy classes in real semisimple Lie groups}},
  journal = {Compositio Mathematica},
  volume = {49},
  number = {3},
  pages = {309--398},
  year = {1983},
  url = {http://www.numdam.org/item?id=CM_1983__49_3_309_0},
  issn = {0010-437X,1570-5846},
  fjournal = {Compositio Mathematica},
  mrclass = {58G25 (22E46 58C27 58G15)},
  mrnumber = {707179},
  mrreviewer = {P.\ C.\ Trombi},
}

@book{Milnor1965,
  author = {Milnor, J.},
  title = {{Lectures on the {$h$}-cobordism theorem}},
  pages = {v+116},
  publisher = {Princeton University Press, Princeton, NJ},
  year = {1965},
  note = {Notes by L. Siebenmann and J. Sondow},
  mrclass = {57.10},
  mrnumber = {190942},
  mrreviewer = {P.\ E.\ Conner},
}

@article{Schwarz1998QuantumCupLength,
  author = {Schwarz, M.},
  title = {{A quantum cup-length estimate for symplectic fixed points}},
  journal = {Inventiones Mathematicae},
  volume = {133},
  number = {2},
  pages = {353--397},
  year = {1998},
  doi = {10.1007/s002220050248},
}

@misc{Gong2025DegenerateFixedPoints,
  author = {Gong, W.},
  title = {{Degenerate symplectic fixed points and Gromov-Witten invariants}},
  howpublished = {arXiv:2507.04191},
  year = {2025},
  doi = {10.48550/arXiv.2507.04191},
  url = {https://arxiv.org/abs/2507.04191},
  primaryclass = {math.SG},
}

@article{RuanTian1995QuantumCohomology,
  author = {Ruan, Y. and Tian, G.},
  title = {{A mathematical theory of quantum cohomology}},
  journal = {Journal of Differential Geometry},
  volume = {42},
  number = {2},
  pages = {259--367},
  year = {1995},
  doi = {10.4310/jdg/1214457234},
}

@incollection{Tamvakis2005IsotropicGrassmannians,
  author = {Tamvakis, H.},
  editor = {Bogomolov, F. and Tschinkel, Y.},
  title = {{Quantum cohomology of isotropic Grassmannians}},
  booktitle = {Geometric Methods in Algebra and Number Theory},
  series = {Progress in Mathematics},
  volume = {235},
  pages = {311--338},
  publisher = {Birkh{\"a}user Boston},
  year = {2005},
  doi = {10.1007/0-8176-4417-2_14},
}

@article{Takens1968,
  author = {Takens, F.},
  title = {{The minimal number of critical points of a function on a compact manifold and the Lusternik-Schnirelman category}},
  journal = {Inventiones Mathematicae},
  volume = {6},
  pages = {197--244},
  year = {1968},
  doi = {10.1007/BF01404825},
  url = {https://doi.org/10.1007/BF01404825},
  issn = {0020-9910,1432-1297},
  fjournal = {Inventiones Mathematicae},
  mrclass = {57.20},
  mrnumber = {236942},
  mrreviewer = {M.\ Klingmann},
}

@article {BanyagaHurtubise2009,
    AUTHOR = {Banyaga, Augustin and Hurtubise, David E.},
     TITLE = {The {M}orse-{B}ott inequalities via a dynamical systems
              approach},
   JOURNAL = {Ergodic Theory Dynam. Systems},
  FJOURNAL = {Ergodic Theory and Dynamical Systems},
    VOLUME = {29},
      YEAR = {2009},
    NUMBER = {6},
     PAGES = {1693--1703},
      ISSN = {0143-3857,1469-4417},
   MRCLASS = {37D15 (37B30 58E05)},
  MRNUMBER = {2563088},
MRREVIEWER = {James\ F.\ Reineck},
       DOI = {10.1017/S0143385708000928},
       URL = {https://doi.org/10.1017/S0143385708000928},
}

@incollection{Giv95,
  author = {Givental, A. B.},
  editor = {Hofer, H. and Taubes, C. H. and Weinstein, A. and Zehnder, E.},
  title = {{A symplectic fixed point theorem for toric manifolds}},
  booktitle = {The Floer Memorial Volume},
  series = {Progress in Mathematics},
  volume = {133},
  pages = {445--481},
  publisher = {Birkh{\"a}user},
  address = {Basel},
  year = {1995},
  doi = {10.1007/978-3-0348-9217-9_19},
  url = {https://math.berkeley.edu/~giventh/papers/tor.pdf},
}
\bibliographystyle{alpha.bst}
\end{document}